\documentclass[]{amsart}
\usepackage[foot]{amsaddr}
\usepackage[a4paper,top=2.5cm,bottom=3cm,inner=3cm,outer=3cm]{geometry}

\usepackage[utf8]{inputenc}

\usepackage{tikz}
\usetikzlibrary{decorations.pathreplacing, decorations.markings, arrows.meta}
\usetikzlibrary{arrows}

\usepackage{amsmath,amsthm}
\usepackage{amsfonts,amssymb,mathrsfs}

\usepackage[bookmarks=false,breaklinks=true]{hyperref}
\usepackage{float}
\usepackage{graphicx}
\usepackage{color}

\usepackage{epigraph}
\definecolor{myurlcolor}{rgb}{0,0,0.7}
\usepackage{hyperref}
\hypersetup{colorlinks,
linkcolor=myurlcolor,
citecolor=myurlcolor,
urlcolor=myurlcolor}
\usepackage[mathscr]{euscript}

\definecolor{darkyellow}{rgb}{0.98,0.98,0}
\definecolor{darkgreen}{rgb}{0,0.9,0}

\definecolor{darkred}{rgb}{0.5,0,0} 

\newcommand{\M}{{\mathcal M}}   

\newcommand{\Z}{{\mathbb Z}}  
\newcommand{\R}{{\mathbb R}}  
\newcommand{\C}{{\mathbb C}}  
\newcommand{\E}{{\mathbb E}}  

\newcommand{\PP}{\mathbb{P}}

\newcommand{\define}[1]{{\bf \boldmath{#1}}}

\theoremstyle{definition}

        \newcommand{\be}{\begin{equation}}
        \newcommand{\ee}{\end{equation}}
        \newcommand{\ba}{\begin{eqnarray}}
        \newcommand{\ea}{\end{eqnarray}}
        \newcommand{\ban}{\begin{eqnarray*}}
        \newcommand{\ean}{\end{eqnarray*}}
        \newcommand{\barr}{\begin{array}}
        \newcommand{\earr}{\end{array}}

\begin{document}
\title{Triangulations of the Sphere}
\author[Baez]{John C.\ Baez} 
\address{Department of Mathematics, University of California, Riverside CA, 92521, USA}
\date{November 5, 2024}
\maketitle

\begin{center}
\includegraphics[width = 20em]{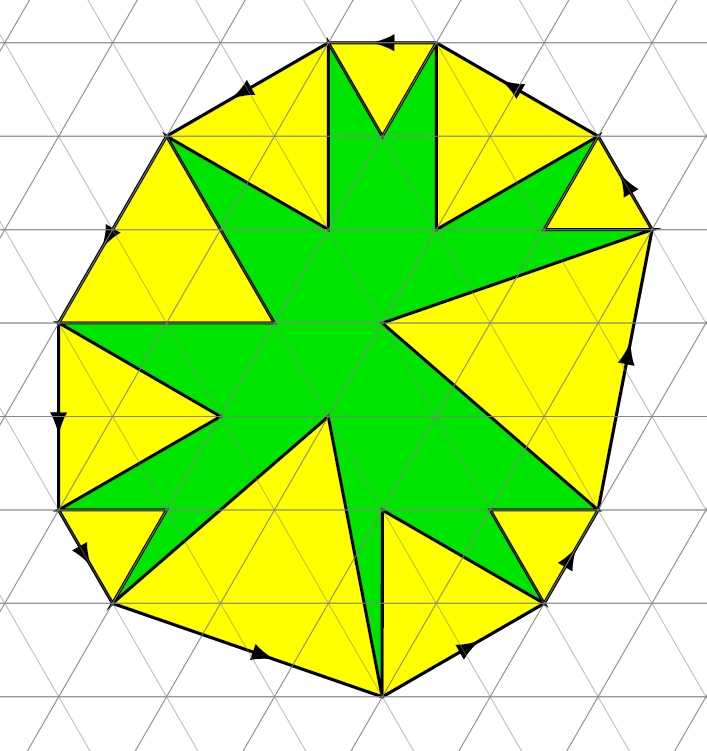} 
\end{center} 

\vskip 1em
If you cut out this awkward-looking green star, you can fold it into a polyhedron where all 11 of its tips meet at one point.  This polyhedron will be tiled by the small gray equilateral triangles, though some may get bent.   There will be 12 points where 5 such triangles meet: the 11 indentations in the star, and the point where all its tips meet.   There are also points where 6 triangles meet.

This may seem like a curiosity, but Bill Thurston discovered some profound facts about stars of this sort \cite{Thurston}.   Schwartz \cite{Schwartz} later tried to clarify Thurston's work, and Engel and Smillie \cite{EngelSmillie} developed it further.  I hope this brief introduction helps more mathematicians discover its charms.

To construct a star of this sort, first draw the lattice of \define{Eisenstein integers} in the complex plane:
\[             \E = \{ a + b \omega \; \vert \; a, b \in \Z  \}  \]
where $\omega = \exp(2 \pi i / 3)$.    Connecting nearest neighbors, you get the small gray equilateral triangles shown above.    Then draw an 11-sided polygon $P$ whose vertices are Eisenstein integers.   Along each edge of $P$, draw a yellow equilateral triangle pointing inward.   Make sure to choose the polygon $P$ so that these yellow triangles touch each other only at its corners.

If you remove these yellow triangles, you are left with an 11-pointed star, like the one shown in green above.  Surprisingly, you can always fold up this star so all its tips meet at one point, forming a convex polyhedron tiled by equilateral triangles: the triangles coming from the lattice of Eisenstein integers.  This polyhedron is homeomorphic to a sphere, so you have constructed a triangulation $T$ of the sphere for which 5 or 6 triangles meet at each vertex.   And Thurston proved a remarkable fact: you can get \emph{all} such triangulations using this method!

``A triangulation of the sphere where 5 or 6 triangles meet at each vertex'' is a purely topological, or combinatorial, concept.  But Thurston noticed that whenever you have such a thing, you can make all its triangles into flat equilateral triangles with the same edge length.  This gives the sphere a geometry.    To be precise, this gives it a flat Riemannian metric except at the points where exactly 5 triangles meet at a vertex.

Your triangulation $T$ thus gives the sphere a flat Riemannian metric, except at finitely many points.  The total angle around each of these points is not $2\pi$, only $5 \pi/3$, because only 5 equilateral triangles meet there.  These points are called `cone points', and we say they have an `angle deficit' of $\pi/3$.   By a discrete version of the Gauss--Bonnet theorem, the total angle deficit must be $4 \pi$, so there must be 12 such cone points.

Thurston showed that up to rescaling,  \emph{any} Riemannian metric on the 2-sphere that is flat except at 12 cone points with angle deficit $\pi/3$ arises from this procedure.  The most famous example is the regular icosahedron, where 5 equilateral triangles meet at each vertex.   Can you see how to draw an 11-sided polygon that gives the regular icosahedron?   If you give up, you can see how Gerard Westendorp did it at the end of this paper.  If you don't believe it, see \cite{Westendorp}.

For deeper results, you can encode your polygon $P$ as a 10-tuple of Eisenstein integers.   If you start at any corner of $P$ and walk all the way around its edges, you get 11 vectors shown as arrows in the figure.   This gives an 11-tuple of Eisenstein integers $(v_1, \dots, v_{11}) \in \E^{11}$.  Since a round trip gets you back where you started, these Eisenstein integers must sum to zero.  Thus, $P$ is determined up to translation by a point $v = (v_1, \dots, v_{10}) \in \E^{10}$.

Since areas depend quadratically on lengths, it is not surprising that there is a real-valued quadratic form $Q$ on $\C^{10}$ such that the number of triangles in the triangulation $T$ equals $Q(v)$.  By a general result in linear algebra, there  is a unique hermitian form $H$ on $\C^{10}$ such that $H(v,v) = Q(v)$.  Engel and Smillie have given a nice explicit formula for this hermitian form \cite[Sec.\ 3]{EngelSmillie}.  It has signature $(1,9)$, meaning that we can find some complex coordinates $z_i$ on $\C^{10}$ such that
\[      Q = |z_1|^2 - |z_2|^2 - \cdots - |z_{10}|^2   .\]
This is reminiscent of 10-dimensional Minkowski spacetime, beloved by string theorists: that is, $\R^{10}$ with the quadratic form $x_1^2 - x_2^2 - \cdots - x_{10}^2$
where the first coordinate describes time and the rest describe space.  Is this just a coincidence?  I do not know.

Thurston studied the so-called `moduli space' $\M$ of all ways of giving a sphere a flat Riemannian metric with 12 cone points with angle deficits of $\pi/3$, where two such metrics count as the same if they are isometric up to rescaling.   He showed $\M$ is open and dense in a larger space 
\[  \overline{\M} = \PP\C^{10}_+/\Gamma. \]
Here
\[         \C^{10}_+ = \{ v \in \C^{10}  \vert \; Q(v) > 0 \}  ,\]
$\PP\C^{10}_+$ is its projectivization (the space of complex lines through the origin in $\C^{10}$ on which $Q$ is positive), and $\Gamma$ is a certain discrete group of linear transformations of $\C^{10}$ preserving both $H$ and the lattice $\E^{10}$.   Since $\C^{10}_+$ has complex dimension $10$, its projectivization $\PP\C^{10}_+$ has complex dimension $9$, and so does $\M$.   The projectivization is what implements the `up to rescaling' in the definition of $\M$.

Thurston also studied the larger space $\overline{\M}$, and showed it is the moduli space of flat Riemannian metrics on the sphere with  \emph{at most} $12$ cone points and angle deficits that are positive integer multiples of $\pi/3$.   So, elements of $\M$ correspond to metrics with 12 distinct cone points, but when two or more cone points collide we get an element of $\overline{\M} - \M$.  $\M$ is a complex manifold, but $\overline{\M}$ is a more singular space, called an `orbifold'.  By Alexandrov's uniqueness theorem \cite{Alexandrov}, any element of $\overline{\M}$ can be realized as a convex polyhedron, unique up to congruence.

Finally, the triangulations themselves correspond not to points of the projectivized space $\overline{\M}$ but to lattice vectors upstairs in $\C^{10}$.  Sitting inside $\C^{10}_+$ we have a discrete set of positive-norm lattice vectors $\C^{10}_+ \cap \E^{10}$, and Thurston showed there is a bijection between oriented convex triangulations of the sphere and the quotient $(\C^{10}_+ \cap \E^{10})/\Gamma$, with the number of triangles given by $Q(v)$ \cite[Thm.\ 0.1]{Thurston}.  These are exactly the triangulations of the sphere where no more than 6 triangles meet at a vertex.    

To see how such triangulations fit into $\overline{\M}$, note that a vertex where $k$ equilateral triangles meet has an angle deficit of $(6-k)\pi/3$. So a vertex where 6 triangles meet is a flat point, a vertex where 5 triangles meet is a cone point of angle deficit $\pi/3$, and a vertex where 4 or 3 triangles meet is a cone point with angle deficit $2\pi/3$ or $\pi$: it is a point where two or three of the original $\pi/3$ cone points have collided. For example, the regular octahedron has six vertices where 4 triangles meet, and the regular tetrahedron has four vertices where 3 triangles meet. Thus, these two Platonic solids do not give points of $\M$, but they do give points of $\overline{\M}$.

\subsubsection*{Acknowledgments} I thank Leo Stein for showing me how to draw the figure and improving the exposition, and Gerard Westendorp for help with the icosahedral case.

\vskip 2em
\begin{center}
\includegraphics[width = 20em]{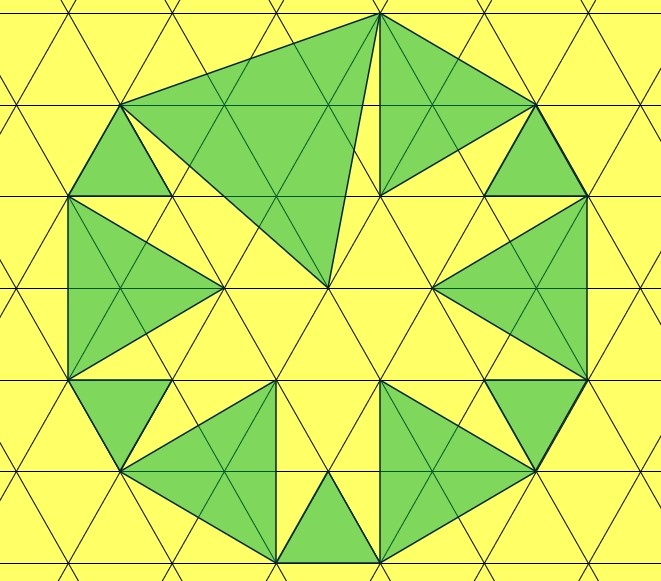} 
\vskip 1em
Gerard Westendorp's construction of the regular icosahedron using Thurston's method.
\end{center} 

\end{document}